\documentclass[nrl,nonblindrev]{informs3}
\usepackage{optidef}
\usepackage{IEEEtrantools}
\usepackage{booktabs}
\usepackage{multirow}
\usepackage{amsmath}
\usepackage{amsfonts}
\usepackage{caption}
\usepackage{graphicx}
\usepackage{subcaption}
\usepackage{algorithm, algpseudocode}
\usepackage{pgfplots}
\pgfplotsset{compat=1.8}
\usepgfplotslibrary{statistics}
\usepackage{pgf} 
\usepackage{pgfplots}
\usetikzlibrary{fit,calc}
\usepgfplotslibrary{external}
\usetikzlibrary{decorations.pathmorphing}
\usetikzlibrary{patterns}
\newcommand{\boxplot}[6]{%
    \filldraw[fill=white,line width=0.2] let \n{boxxl}={#1-0.1}, \n{boxxr}={#1+0.1} in (axis cs:\n{boxxl},#3) rectangle (axis cs:\n{boxxr},#4);   
    \draw[line width=0.2mm, color=red] let \n{boxxl}={#1-0.1}, \n{boxxr}={#1+0.1} in (axis cs:\n{boxxl},#2) -- (axis cs:\n{boxxr},#2);             
    \draw[line width=0.2mm] (axis cs:#1,#4) -- (axis cs:#1,#6);                                                                           
    \draw[line width=0.2mm] let \n{whiskerl}={#1-0.025}, \n{whiskerr}={#1+0.025} in (axis cs:\n{whiskerl},#6) -- (axis cs:\n{whiskerr},#6);        
    \draw[line width=0.2mm] (axis cs:#1,#3) -- (axis cs:#1,#5);                                                                           
    \draw[line width=0.2mm] let \n{whiskerl}={#1-0.025}, \n{whiskerr}={#1+0.025} in (axis cs:\n{whiskerl},#5) -- (axis cs:\n{whiskerr},#5);        
}
\usepackage{rotating}
\DoubleSpacedXI 

\usepackage{endnotes}
\let\footnote=\endnote

\usepackage{natbib}
 \bibpunct[, ]{(}{)}{,}{a}{}{,}%
 \def\bibfont{\small}%
 \def\BIBand{and}%

\TheoremsNumberedThrough     
\ECRepeatTheorems

\EquationsNumberedThrough    

\MANUSCRIPTNO{xxxx} 

\begin{document}


\RUNAUTHOR{Kiani, Eksioglu, Isik, Thomas, Gilpin}

\RUNTITLE{Appointment Postponement in Scheduling Patients}

\TITLE{Evaluating Appointment Postponement in Scheduling Patients at a Diagnostic Clinic}

\ARTICLEAUTHORS{%
\AUTHOR{Mahsa Kiani, Burak Eksioglu, Tugce Isik, Alexandria Thomas}
\AFF{Department of Industrial Engineering, Clemson University, Clemson, SC 29634}
\AUTHOR{John Gilpin}
\AFF{Radiology Department, Hillcrest Memorial Hospital, Simpsonville, SC 29681}
\AFF{\EMAIL{mkiani@clemson.edu}, \EMAIL{burak@clemson.edu}, \EMAIL{tisik@clemson.edu}, \EMAIL{amt3@clemson.edu}, \EMAIL{jgilpin@ghs.org}}

} 

\ABSTRACT{Diagnostic clinics are among healthcare facilities that suffer from long waiting times which can cause medical issues and lead to increases in patient no-shows. Reducing waiting times without significant capital investments is a challenging task. We tackle this challenge by proposing a new appointment scheduling policy for such clinics that does not require significant investments. The clinic in our study serves outpatients, inpatients, and emergency patients. Emergency patients must be seen on arrival, and inpatients must be given next day appointments. Outpatients, however, can be given later appointments. The proposed policy takes advantage of this by allowing the postponement of the acceptance of appointment requests from outpatients. The appointment scheduling process is modeled as a two-stage stochastic programming problem where a portion of the clinic capacity is allocated to inpatients and emergency patients in the first stage. In the second stage, outpatients are scheduled based on their priority classes. After a detailed analysis of the solutions obtained from the two-stage stochastic model, we develop a simple, non-anticipative policy for patient scheduling. We evaluate the performance of this proposed, easy-to-implement policy in a simulation study which shows significant improvements in outpatient indirect waiting times.}



\maketitle

%
\section{Introduction}
\label{section:intro}
In today's healthcare systems, the increasing demand for appointments combined with a shortage of physicians has led to challenges for healthcare providers to give timely appointments to patients. To achieve good medical outcomes, offering timely appointments is important \citep{gupta}. \cite{gupta} classify waiting time of patients into two categories. They define \textit{direct waiting time} as the time the patient waits in the healthcare facility on the day of appointment and \textit{indirect waiting time} as the time between the day the patient requests an appointment and the appointment day. Unfortunately, long indirect waiting times are common in practice. For instance, \cite{kesling} reported that 84$\%$ of patients in Veterans Affairs (VA) hospitals wait more than 14 days to see a physician. In addition to the medical issues that long indirect waiting times cause, they can also lead to increases in patient no-shows \citep{green} which have significant effect on annual revenues \citep{moore}. Thus, healthcare managers face the challenge of improving their appointment systems to decrease waiting times and no-shows without incurring major capital costs.

Diagnostic clinics are among the healthcare facilities that generally suffer from long indirect waiting times \citep{mccarthy}. One such clinic is the Radiology Department at Prisma Health, our collaborator on this study. The clinic provides service to outpatients, inpatients, and emergency patients. The requests for appointments are handled on a first-come-first-served (FCFS) basis. The emergency patients are the highest priority group, followed by inpatients and then outpatients. The outpatients are further categorized into a number of priority classes based on co-morbidities and chronic conditions. The emergency patients are seen as soon as they arrive if there is capacity or immediately referred to another clinic. The inpatients are either given a next day appointment during regular hours or seen during overtime hours. The clinic prefers to offer appointments to outpatients within a few days. However, under the current system, the average indirect waiting time for outpatients is about one week. \cite{luo} provide other examples where indirect waiting times for outpatients are negatively impacted by the arrival of higher priority inpatients and emergency patients. A possible strategy to reduce the indirect waiting times for outpatients is to allocate a part of the overall capacity for emergency patients to dampen their impact on the overall system. Similarly, a portion of the capacity can also be reserved for inpatients. However, this strategy can result in unused capacity. Meanwhile, the limited available capacity may not allow providers to serve some of the more urgent outpatients in an acceptable time period. Thus, finding ways to utilize the unused portion of the capacity reserved for inpatients and allocating just enough capacity for emergency patients are important problems.

The clinic currently makes all acceptance and referral decisions upon the arrival of appointment requests. This causes some high priority outpatients to be referred to other clinics while some of the capacity reserved for inpatients goes unused. As a solution, we propose postponing the acceptance of outpatient requests. In other words, the decision regarding acceptance or referral of an outpatient is not taken upon arrival of an appointment request but is revisited after the inpatient schedules are realized. This postponement will enable the scheduling of higher priority arrivals sooner and also allow for better utilization of the unused capacity reserved for inpatients. Note that postponement does not allow one to utilize the potential unused capacity allocated for emergency patients, because we do not have the one day buffer which is the case for inpatients. Thus, it is critical to allocate the right amount of capacity for emergency patients. 

The majority of the outpatients prefer to get an immediate response from the clinic regarding their appointment request. However, the clinic is willing to keep outpatient appointment requests in an \textit{acceptance queue} for a reasonable amount of time. While some patients may leave for an alternative healthcare facility, the clinic believes that most of the outpatients will be amenable to waiting in the acceptance queue if it means their total indirect waiting time will be shorter. Still, the clinic is not open to keeping the outpatients in the acceptance queue more than 72 hours. 

To that end, we develop a two-stage, postponable acceptance appointment model which first allocates the total regular-time capacity among different groups of patients and then schedules appointments. Outpatient appointment requests are either scheduled during regular hours or referred to another clinic. The objective is to minimize the expected total cost over the planning horizon. The remainder of the paper is organized into five additional sections. Section \ref{section:literature} provides a review of the relevant literature. In Section \ref{section:formulation} the problem is formally defined and a notation is provided along with a two-stage stochastic programming (TSSP) formulation. Section \ref{section:Solution Approach} explains how the problem is solved. Specifically, the details of our sample average approximation (SAA) and decomposition-based branch-and-bound (DBB) algorithm are provided. Section \ref{section:results} shows the results of our extensive experiments and sensitivity analysis. Finally, Section \ref{section:conclusion} concludes the paper with some managerial insights, highlights some of the limitations, and provides directions for future research.

\section{Literature Review}\label{section:literature}
Our study is related to four streams of literature. In the following paragraphs we provide brief reviews of the related literature on ($i$) patient scheduling, ($ii$) acceptance postponement, ($iii$) solution approaches for TSSP, and ($iv$) revenue management. We highlight how our study differs from those in the literature and summarize our contributions. 

The scheduling of patients with different priority classes and medical resource allocation to these classes has gotten a lot of attention in recent years, as evidenced by the large number of papers in the literature \citep{patrick, qu2013two, berg2014optimal, feldman2014appointment, kong2015appointment, jiang2017integer}. \cite{ahmadi} provide a comprehensive review of recent analytical and numerical studies in the area of outpatient scheduling. Some of these studies consider inpatients and emergency patients in addition to outpatients, where the arrival of inpatients and emergency patients are modeled as random events that interrupt the system \citep{patrick2007improving, erdogan2013dynamic, erdogan2015online}. \cite{deglise2018capacity} provide a capacity allocation plan to  minimize the indirect waiting time of higher priority patients across an integrated network of care services. On the other hand, scheduling of outpatients in the presence of emergency and inpatient arrivals is studied via appointment scheduling, but not capacity planning, in diagnostic clinics by \cite{green,sickinger}, and \cite{bhattacharjee}.
\cite{green} discuss scheduling of patients in a diagnostic clinic where a certain number of outpatients are already scheduled. They assume that emergency patients arrive randomly throughout the day, and they have to be seen as soon as they arrive. They specify which patient to schedule next when both outpatients and inpatients are waiting for appointments. 
\cite{bhattacharjee} study outpatient and inpatient scheduling problems with non-homogeneous mean service times considering punctuality and no-show rates. Reserving a part of capacity for emergency arrivals or inpatients is shown to be beneficial to decrease the waiting time of urgent patients \citep{patrick2007improving}. \cite{tang} apply a robust optimization model in a surgery department to decide how much capacity to allocate for elective surgeries and emergency surgeries when the demand is uncertain.
The decision regarding acceptance or rejection of patients depends on their priorities and available capacity. \cite{akhavizadeganl} develop a finite-horizon Markov decision process to schedule appointments considering choice behavior and no-show rate of patients. Patients provide their preference for a specific physician and time of appointment. The decision to accept or reject walk-in patients is based on already scheduled patients who called-in. The main difference between our study and those mentioned above is the timing of the decisions regarding acceptance, rejection, or referral of outpatient appointment requests. 

The concept of acceptance postponement is developed and discussed in some manufacturing settings but not so much in service settings. For example, \cite{kang} present a model for a manufacturing system with postponable acceptance and assignment in make-to-order settings, where postponement is applied to both acceptance and assignment. In their model, acceptance of some orders may be deferred or cancelled to wait for more profitable orders. They show that by applying this model the total profit of the system improves. In a study by \cite{gao2012managing}, some low-priority orders are rejected or the acceptance decision is postponed to reserve inventory for higher priority orders. \cite{bassamboo} provide one of the few studies of applying postponement in a service system. They consider a call center routing problem that assigns arrivals right after acceptance or after some waiting period. However, acceptance of calls have to be made at the time of arrival. Moreover, acceptance and assignment decisions are made at the same time if there is an available agent. 
The two closest studies to ours are by \cite{balasubramanian} and \cite{patrick}.
\cite{balasubramanian} consider both open-access and prescheduled appointments in their settings. They compute how much of a physician's workload should be allocated to prescheduled appointments. However, scheduling of patients occurs upon their arrivals.
In contrast, the study by \cite{patrick} considers the acceptance of some of the requests to be postponed. They consider a dynamic system which schedules multiple priority classes of outpatients with the goal of decreasing indirect waiting times when the daily outpatient capacity is fixed. In their model, once the acceptance decisions are made, the remaining requests are deferred to the next day and may be accepted later. However, they do not keep track of the number of days that the decisions are deferred. We postpone the acceptance and scheduling of outpatients in our setting as well. However, our study considers the following concepts that are not included in \cite{patrick} study.
First, we consider the cost of postponing the acceptance decisions, which depends on the amount of time outpatients wait in the acceptance queue and their priority classes. Second, we consider an abandonment probability which relies on the outpatient's priority class and the amount of time s/he has waited in the acceptance queue. Third, we consider how the postponing of outpatients affects the capacity allocation and scheduling of inpatients and emergency patients.
To the best of our knowledge, our study is the first one that introduces a capacity allocation and postponement model for patient scheduling. 

As discussed in Section \ref{section:formulation}, we formulate our problem as a TSSP and develop SAA and DBB approaches, as detailed in Section \ref{section:Solution Approach}, to solve the problem. A well-known approach to solve TSSP is stochastic Benders decomposition which is also known as the L-shaped method. \cite { van1969shaped} were the first ones to use the L-shaped method to solve TSSP problems. In their formulation, the first and second stage variables were all continuous. \cite{laporte1993integer} allowed integer first and/or second stage variables in their setting by incorporating a branch-and-bound procedure. \cite{ahmed2004finite}  proposed the DBB algorithm by branching on tender variables that are the product of first stage variables with the technology matrix for problems with integer variables in the second stage. In our study, we first replace the original objective function by a SAA function and then apply the DBB algorithm to be able to solve realistic size problems.

While our study does not directly contribute to the revenue management literature, there are similarities. Revenue management is defined as the management of perishable assets \citep{weatherford1992taxonomy}. Examples of perishable assets include hotel rooms, rental cars, and airplane seats. Revenue management of these perishable assets includes the process of allocating a fixed capacity to the right customer at the right time at the right price \citep{smith1992yield}. One of the studies which is close to ours is where they allocate the scarce inventory to stochastic demand for multiple fare classes so as to maximize the total expected revenue \citep{bertsimas2003revenue}. The structure of optimal policy is estimated by solving an approximate dynamic
program. Revenue management decisions are made upon arrivals but considering anticipated future requests. In this perspective, our study is different since we consider the possibility of postponing the decisions. 

\section{Problem Definition and Formulation} \label{section:formulation}
As mentioned in Section \ref{section:intro}, the diagnostic clinic in our study receives appointment requests from outpatients, inpatients, and emergency patients. Currently, almost all of the outpatient appointment requests are accepted or referred to another clinic as soon as the request arrives. The only exception to this are those requests that are received via fax which constitute a small fraction of all requests. The clinic responds to the faxed requests by the end of the business day. We, on the other hand, develop a TSSP that allows the postponement of all outpatient requests.

Outpatients are categorized into $J$ priority classes ($j=1,...J$). Parts of the regular-time capacity are allocated for inpatients and emergency patients. The capacity reserved for inpatients can be used for outpatients only if it is unused after inpatients are scheduled. Emergency patients that arrive throughout the day are either seen upon arrival or immediately referred to another clinic. Inpatients that arrive throughout the day are either given a next day appointment during regular hours upon arrival or seen during overtime hours. Outpatient requests that arrive each day are kept in the acceptance queue. In other words, the acceptance and scheduling decisions of lower priority outpatients can be postponed while waiting for inpatients, emergency patients, or higher priority outpatients. To facilitate the formulation of our model we define the parameters shown in Table \ref{table:parameters} and the variables shown in Table \ref{table:variables}.
\begin{table}[h!]
\footnotesize
\centering
\begin{tabular}{rcl}
 \hline
 \multicolumn{3}{l}{Parameters} \\
 \hline
 $T$ & : & length of the planning horizon ($t=1,2,...,T$) \\
 $T^a$ & : & length of the booking horizon ($t^a=1,2,...,T^a$) \\
 $T^u$ & : & maximum number of days an outpatient waits in the acceptance queue ($t^u=1,2,...,T^u$) \\
 $K$ & : & daily regular-time capacity of the clinic \\
 $p_{jt^u}$ & : & proportion of type $j$ outpatients who stay in the acceptance queue one more day after \\
 & & having waited for ($t^u$-1) days \\
 $a_{jt^u}$ & : & cost of a type $j$ outpatient leaving the acceptance queue after waiting for $t^u$ days \\
 $b_{jt^ut^a}$ & : & cost of giving an appointment to a type $j$ outpatient $t^a$ days later when the patient\\
 & & has waited for $t^u$ days in the acceptance queue \\
 $c^O_{jt^u}$ & : & cost of referring a type $j$ outpatient to another clinic when the patient has waited \\
 & & for $t^u$ days in the acceptance queue \\
 $c^I$ & : & cost of seeing an inpatient during overtime hours \\
 $c^E$ & : & cost of referring an emergency patient to another clinic \\
\hline
\end{tabular}
\caption{Problem Parameters}
\label{table:parameters}
\end{table}

In our proposed system, a scheduler observes the number of inpatients ($D^I_{t}$) and outpatients ($D^O_{jt}$) that have arrived during the day and the available capacity in each future day of the booking horizon. If the daily inpatient arrivals exceed the allocated capacity ($K\alpha^I$), they are handled during overtime hours which incurs additional cost ($c^I$). If any of the capacity allocated to inpatients is not used then it can be allocated to an outpatient from the acceptance queue. However, the capacity reserved for emergency patients ($K\alpha^E$) is never used for inpatients or outpatients. If an emergency patient arrives when the allocated capacity is full then s/he is immediately referred to another clinic. An outpatient who has been in the acceptance queue for $T^u$ days is referred to another clinic.
\begin{table}[h!]
\footnotesize
\centering
\begin{tabular}{rcl}
 \hline
 \multicolumn{3}{l}{Random Variables} \\
 \hline
 $D^O_{jt}$ & : & number of type $j$ outpatients that arrive at the clinic during day $t$ \\
 $D^I_{t}$ & : & number of inpatients that arrive at the clinic during day $t$ \\
 $D^E_{t}$ & : & number of emergency patients that arrive at the clinic during day $t$ \\
 \hline
 \multicolumn{3}{l}{Decision Variables} \\
 \hline
 $\alpha^I$ & : & percentage of total capacity $K$ reserved for inpatients \\
 $\alpha^E$ & : & percentage of total capacity $K$ reserved for emergency patients \\
 $y^{O}_{jtt^ut^a}$ & : & number of type $j$ outpatients who are given an appointment in day $t$ for $t^a$ days later \\
 & & after waiting for $t^u$ days in the acceptance queue ($t^u \leq t$) \\
 $r_{jtt^u}$ & : & number of type $j$ outpatients who are referred to another clinic in day $t$ after \\
 & & waiting for $t^u$ days in the acceptance queue ($t^u \leq t$) \\
 $Q_{jtt^u}$ & : & number of type $j$ outpatients in day $t$ who have been waiting in the acceptance  \\
 & & queue for $t^u$ days ($t^u \leq t$) \\
 $K^O_{tt^a}$ & : & available capacity for outpatients $t^a$ days after day $t$ \\
 \hline
\end{tabular}
\caption{Problem Variables}
\label{table:variables}
\end{table}

Based on analysis of historical data and our conversations with Prisma Health, patient arrivals are independent Poisson processes. Thus, we model $D^O_{jt}, D^I_{t}$ and $D^E_{t}$ as truncated Poisson distributions with rates $\lambda_j, \lambda^I$ and $\lambda^E$, respectively. The evolution of $Q_{jtt^u}$, the number of outpatients in the acceptance queue, is captured by the following equations: 
{\footnotesize
\begin{IEEEeqnarray}{rl"l}
\IEEEyesnumber\label{eq:QueueSize}
\IEEEyessubnumber*
Q_{jt1} &= p_{j1}D^O_{jt} - \sum_{t^a=1}^{T^a}{y^O_{jt1t^a}} - r_{jt1}, &\forall j,t, \label{eqn:Q1} \\
Q_{j(t+1)(t^u+1)} &= p_{j(t^u+1)}Q_{jtt^u} - \sum_{t^a=1}^{T^a}{y^O_{i(t+1)(t^u+1)t^a}} - r_{j(t+1)(t^u+1)}, &\forall j,t,t^u, t\neq T, t^u \neq T^u, t^u \leq t, \label{eqn:Qt}
\end{IEEEeqnarray}
}
\noindent where equation (\ref{eqn:Q1}) is a special case of equation (\ref{eqn:Qt}) with $t^u = 1$. These two equations simply state that the number of outpatients in the next day will be equal to the number of outpatients who are not scheduled or referred yet and remained in the queue for one more day.  

We also need to maintain an accurate account of the remaining regular-time capacity. This can be achieved by the following equations where (\ref{eqn:K1}) is for the beginning of the planning horizon, (\ref{eqn:KTa}) for the end of the booking horizon, and equation (\ref{eqn:Kt}) for other days during the planning and booking horizons. At the beginning of the planning horizon and the end of the booking horizon we have full capacity for outpatients since no body is scheduled in these days yet. In the remaining days, the available capacity on day ($t+1$) is available capacity of day $t$ minus the scheduled appointments for that day.  
{\footnotesize
\begin{IEEEeqnarray}{rl"l}
\IEEEyesnumber\label{eq:capacity}
\IEEEyessubnumber*
K^O_{1t^a} &= K(1-\alpha^I - \alpha^E), &\forall{t^a}, \label{eqn:K1} \\
K^O_{tT^a} &= K(1-\alpha^I - \alpha^E), &\forall{t}, \label{eqn:KTa} \\
K^O_{(t+1)t^a} &= K^O_{t(t^a+1)} - \sum_{j=1}^{J}{\sum_{t^u=1}^{T^u}{y^O_{jtt^u(t^a+1)}}}, &\forall{t,t^a}, t \neq T, t^a \neq T^a. \IEEEyessubnumber \label{eqn:Kt}
\end{IEEEeqnarray}
}
The postponable acceptance appointment system can now be formulated as the following TSSP, named (2SIP). Since capacity allocations have to be made prior to the realization of patient arrivals, $\boldsymbol{\alpha}=(\alpha^I,\alpha^E)$ are the first-stage decision variables. On the other hand, the appointments depend on patient arrivals. Thus, $\boldsymbol{y^O, r, Q,}$ and $\boldsymbol{K^O}$ are the second-stage variables. 
{\footnotesize
\begin{mini!}|s|[0]
{\boldsymbol{\alpha}}
{\mathbb{E_{\boldsymbol{\omega} \in \boldsymbol{\Omega}}}[C(\boldsymbol{\alpha},\boldsymbol{\omega})]\label{S1:obj}}
{\label{stage1}}
{\text{(2SIP)} \quad C^*=}
\addConstraint{\alpha^I + \alpha^E}{\leq 1 \label{S1:c1}}{}
\addConstraint{\alpha^I , \alpha^E}{\geq 0. \label{S1:c2}}{}
\end{mini!}
}
The model minimizes the expected total cost associated with appointment scheduling. Note that $\boldsymbol{\omega}=\{ (D^O_{1t},...,D^O_{Jt},D^I_t,D^E_t)$ for $t=1,...,T \}$ is a joint scenario for the planning horizon. We assume that there is no cost for capacity allocation. The objective of the second stage is to minimize the cost associated with scheduling patient appointments. As shown in Table \ref{table:parameters}, costs are incurred when outpatients abandon the acceptance queue, outpatients are given late appointments, outpatients are referred to another clinic, inpatients are seen during overtime hours, and emergency patients are referred to another clinic. Thus, $C(\boldsymbol{\alpha},\boldsymbol{\omega})$ is the objective function value of the second-stage problem given below:
{\footnotesize
\begin{mini!}|s|[0]
{\boldsymbol{y^O},\boldsymbol{r},\boldsymbol{Q},\boldsymbol{K}^O}
{\Bigg\{\sum_{t=1}^{T} \Big(\sum_{j=1}^{J}{\sum_{t^u=1}^{T^u}{\sum_{t^a=1}^{T^a}{b_{jt^ut^a} y_{jtt^ut^a}^O}}} +\sum_{j=1}^{J}{\sum_{t^u=1}^{T^u}{{c^{O}_{jt^u}r_{jtt^u}}}} + a_{j{t^u}}(1 - p_{jt^u})Q{jtt^u}  \label{S2:obj}}
{\label{stage2}}
{C(\boldsymbol{\alpha},\boldsymbol{\omega})=}
\addConstraint{\nonumber + c^{I}(1-z^I_t) (D^I_t - \alpha^I K) + c^{E}(1-z^E_t)(D^E_t - \alpha^E K)\Big)\Bigg\}}
\addConstraint{(\ref{eqn:Q1})-(\ref{eqn:Kt})}{\label{S2:cons0}}{}
\addConstraint{z^{I}_t K}{\geq \alpha^I K - D^{I}_t, \quad \label{S2:cons1}}{\forall t}
\addConstraint{z^{I}_t D^{I}_t}{\leq \alpha^I K, \quad \label{S2:cons2}}{\forall t}
\addConstraint{z^{E}_t K}{\geq \alpha^E K - D^{E}_t, \quad \label{S2:cons3}}{\forall t}
\addConstraint{z^{E}_t D^{E}_t}{\leq \alpha^E K, \quad \label{S2:cons4}}{\forall t}
\addConstraint{\sum_{j=1}^{J}{\sum_{t^u=1}^{T^u}{{y^{O}_{jtt^u1}}}} - z^I_t(\alpha^I K - D^I_t)}{\leq K^{O}_{t1}, \quad \label{S2:cons5}}{\forall t}
\addConstraint{\sum_{j=1}^{J}{\sum_{t^u=1}^{T^u}{y^{O}_{jtt^ut^a}}}}{\leq K^{O}_{tt^a}, \quad \label{S2:cons6}}{\forall t, t^a=2,...,T^a}
\addConstraint{Q_{jtT^u}}{= r_{jtT^u}, \quad \label{S2:cons7}}{\forall j,t, t \geq T^u}
\addConstraint{z^{I}_t, z^{E}_t}{\in \{0,1\}, \quad \label{S2:cons8}}{\forall{t}}
\addConstraint{y^{O}_{jtt^ut^a}, r_{jtt^u}, Q_{jtt^u}, K^O_{tt^a}}{\in \mathbb{Z^+}, \quad \label{S2:cons9}}{\forall j,t,t^a, t^u \leq t.}
\end{mini!}
}
To model whether or not demand exceeds capacity, we introduce binary variables $z^{I}_t$ and $z^{E}_t$. We let $z^{I}_t=1$ if $D^{I}_t \leq \alpha^I K$ at time $t$ and $z^{I}_t=0$ otherwise. Similarly, $z^{E}_t=1$ if $D^{E}_t\leq \alpha^E K$ and 0 otherwise. Constraints (\ref{S2:cons1})-(\ref{S2:cons4}) ensure that $z^{I}_t$ and $z^{E}_t$ take on the correct values depending on whether or not demand is less than the corresponding capacity. Constraint set (\ref{S2:cons5}) ensures that the total number of next day appointments given to outpatients does not exceed the remaining capacity for outpatients plus the unused capacity that was reserved for inpatients. Constraint set (\ref{S2:cons6}) is similar to (\ref{S2:cons5}), \textit{i.e.}, it ensures that the total number of outpatient appointments does not exceed the remaining capacity on the subsequent days. The only difference is that in (\ref{S2:cons5}) we also have the unused capacity that was initially allocated for inpatients which can now be used for outpatients. Constraint set (\ref{S2:cons7}) ensures that patients do not wait more than $T^u$ days in the queue. Finally, constraints (\ref{S2:cons8}) and (\ref{S2:cons9}) are the binary and integrality constraints.

Note that, when solving the first-stage problem (2SIP), the objective function (\ref{S2:obj}) and the constraint set (\ref{S2:cons5}) are  nonlinear. However, we will not linearize these since our approximation and decomposition approaches will not require solving (2SIP) directly. Instead, we will reformulate the problem as described in Section \ref{section:Solution Approach}. 

\noindent \textbf{Limitations of the model:} One of the limitations of our TSSP model is that is anticipative, \textit{i.e.}, it relies on knowing the demand for the whole planning horizon. Another limitation is that the model assumes the system is initially empty. Also, the model is considering a finite planning horizon which may lead to end-of-horizon effects. We address these limitations to some extend as discussed later in Sections \ref{section:Solution Approach} and \ref{section:results}.

\section{Solution Approach}\label{section:Solution Approach}
Due to the curse of dimensionality, solving (2SIP) as presented in Section \ref{section:formulation} is impractical. To overcome this complexity, we develop a sample average approximation (SAA) approach to generate tight upper and lower bounds. The SAA procedure generates a random sample $\boldsymbol{\omega}^1, \boldsymbol{\omega}^2, \ldots, \boldsymbol{\omega}^S$ of $S$ scenarios from $\Omega$, the set of all possible scenarios, and solves $M$ replications of the following deterministic SAA problem: 
{\footnotesize
\begin{mini!}|s|[0]
{\boldsymbol{\alpha}}
{\frac 1 S \sum_{s=1}^{S}{C(\boldsymbol{\alpha},\boldsymbol{\omega}^s)} \label{SAA:obj}}
{\label{app2slp}}
{\hat{\text{(2SIP)}} \quad \hat{C}_S=}
\addConstraint{(\ref{S1:c1}), (\ref{S1:c2}).}{}{}
\end{mini!}
}
Note that $\hat{C}_S \rightarrow C^*$ as $S \rightarrow \infty$, and estimates of the optimal first-stage solutions for the original stochastic problem can be obtained by solving this deterministic version \citep{verweij}. Algorithm \ref{al:SAA} below formalizes our proposed SAA approach. As shown in the algorithm, the average of the $M$ replications ($\bar{C}_S$) provides a statistical lower bound for $C^*$. For each solution to $\hat{\text{(2SIP)}}$ from the $M$ replications, the second-stage problem (\ref{stage2}) is solved using a larger sample size $S'$. Among this larger sample, the one with the smallest objective value ($\hat{C}_{S'}(\boldsymbol{\hat{\alpha}}^*)$) is our statistical upper bound for $C^*$. We also calculate the variances of the lower and upper bound estimates, \textit{i.e.}, $\sigma^2_{\bar{C}_{S}}$ and $\sigma^2_{\hat{C}_{S'}(\boldsymbol{\hat{\alpha}}^*)}$, respectively. The proofs of the estimation of lower and upper bounds are provided by \cite{mak} and \cite{verweij}, and thus, omitted here. The algorithm increases the sample sizes $S$ and $S'$ until the optimality gap and the variance of the gap estimator are small.
\begin{algorithm}[h!]
\footnotesize
\caption{Sample Average Approximation (SAA)}\label{al:SAA}
\begin{algorithmic}\vspace{0.2cm}
\State \textbf{Step 1}: Initialize $S$, $S'$, and $M$; 
\State \textbf{Step 2}: \textbf{For}  $m = 1,..., M$
    \State \hspace{0.2in} Solve $\hat{\text{(2SIP)}}$ using DBB to obtain objective values $\hat{C}^m_S$ and solutions $\boldsymbol{\hat{\alpha}}^m$;
\State \textbf{Step 3}: Calculate $\bar{C}_S = \frac{1}{M} \sum_{m=1}^{M} \hat{C}^m_S $ and $\sigma^2_{\bar{C}_{S}} = \frac{1}{M(M-1)} \sum_{m=1}^{M}(\hat{C}^m_S - \bar{C}_S)^2$;
\State \textbf{Step 4}: \textbf{For}  each $\boldsymbol{\hat{\alpha}}^m$
    \State \hspace{0.2in} Solve (\ref{stage2}) and compute $\hat{C}_{S'} = \frac{1}{S'} \sum_{s=1}^{S'}{C(\boldsymbol{\hat{\alpha}}^m,\boldsymbol{\omega}^s)}$ and $\sigma^2_{\hat{C}_{S'}(\boldsymbol{\hat{\alpha}})} = \frac{1}{S'(S'-1)} \sum_{s=1}^{S'}(C(\boldsymbol{\hat{\alpha}}^m,\boldsymbol{\omega}^s) - \hat{C}_{S'})^2$;
\State \textbf{Step 5}: Let $\boldsymbol{\hat{\alpha}}^* = \arg \min \Big\{ \hat{C}_{S'}(\boldsymbol{\hat{\alpha}}) : \boldsymbol{\hat{\alpha}} \in \{\boldsymbol{\hat{\alpha}}^1, \ldots, \boldsymbol{\hat{\alpha}}^M \} \Big\}$;
\State \textbf{Step 6}: Calculate $\Delta_C = \hat{C}_{S'}(\boldsymbol{\hat{\alpha}}^*)- \bar{C}_S$ and $\sigma^2=\sigma^2_{\bar{C}_{S}}+\sigma^2_{\hat{C}_{S'}(\boldsymbol{\hat{\alpha}}^*)}$;
\State \textbf{Step 7}: \textbf{If} ($\Delta_C < \epsilon$ and $\sigma^2 < \epsilon$) \textbf{then} report $\boldsymbol{\hat{\alpha}}^*$ as the optimal solution and terminate;
\State \hspace{0.2in} \textbf{Else} increase $S$ and $S'$ and go back to Step 2.
\end{algorithmic}
\end{algorithm}

Solving $\hat{\text{(2SIP)}}$ in Algorithm \ref{al:SAA}, while easier than solving (2SIP), is still a challenging task for large $S$. To that end, we developed a decomposition based branch-and-bound (DBB) algorithm which was originally proposed by \cite{ahmed2004finite} to solve TSSP models with continuous first-stage and discrete second-stage variables. The main idea behind DBB is to partition the search space to efficiently identify candidate solutions \citep{ahmed2002sample,ahmed2004finite}. To be able to implement DBB and ensure convergence, the following assumptions must be satisfied (all of which are satisfied for $\hat{\text{(2SIP)}}$): (A1) The uncertain parameter $\omega$ follows a discrete distribution with finite support. (A2) The first-stage constraint set is nonempty and compact. (A3) The second-stage variables are purely integer. (A4) The technology matrix is deterministic. (A5) For each scenario the second-stage problem is bounded. (A6) For each scenario, the second-stage constraint matrix is integral. We reformulate $\hat{\text{(2SIP)}}$ as follows: 
{\footnotesize
\begin{mini!}|s|[0]
{\boldsymbol{\chi}}
{f(\boldsymbol{\chi}) \label{TPobj}}
{\label{reform2sip}}
{\text{(TP)} \quad }
\addConstraint{\boldsymbol{\chi} \in X,}{\label{TP:c1}}{}
\end{mini!}
}
where $f(\boldsymbol{\chi}) = \frac 1 S \sum_{s=1}^{S}{\Psi^s(\boldsymbol{\chi})}$, $\Psi^s(\boldsymbol{\chi}) = \min \{ f^s \boldsymbol{y} \mid D^s \boldsymbol{y} \geq h^s + \boldsymbol{\chi}, \boldsymbol{y} \in Y\cap \mathbb{Z} \}$, and $X = \{ \boldsymbol{\chi} \mid \boldsymbol{\chi} =T\boldsymbol{\alpha}, \text{ with } (\ref{S1:c1}) \text{ and } (\ref{S1:c2}) \}$. In the stochastic programming literature, the matrix $T$ is known as the \textit{technology matrix} and variables $\boldsymbol{\chi}$ as the \textit{tender variables} that link the first- and second-stage problems. Note that for our problem $T$ is deterministic, \textit{i.e.}, it is independent of the scenario observed. The term $\Psi^s(\boldsymbol{\chi})$ is essentially a compact representation of the second-stage problem given by formulation (\ref{stage2}) where $\boldsymbol{y}$ represents the collection of all second-stage decision variables (\textit{i.e.}, $\boldsymbol{y} = (\boldsymbol{y^O, r, Q, K^O})$), $f^s$ represents the objective function (\ref{S2:obj}), and $D^s$, $h^s$, and $Y$ represent the constraints (\ref{S2:cons0})-(\ref{S2:cons9}) with $D^s$ corresponding to the scenario dependent coefficients, $h^s$ the scenario dependent constants, $Y$ the scenario independent constraints, and $T$ the scenario independent parts of the constraint set which include the first-stage variables. This reformulation allows us to consider a larger number of scenarios in Step 2 of Algorithm \ref{al:SAA}. More specifically, the DBB algorithm below enables us to avoid solving $\hat{\text{(2SIP)}}$ directly. Instead of the first-stage variables, we search the space of the tender variables for global optima. The search space of $\chi$ is partitioned into subsets of the form $\prod_j(l_j,u_j]$, for each component $j$ of $\chi$ where $l_j$ is a point at which the second-stage value function ($\Psi^s(\cdot)$) may be discontinuous \citep{ahmed2004finite}. By branching this way, we isolate subsets over which the second-stage value function is constant. Thus, we can solve $\hat{\text{(2SIP)}}$ exactly. 

\begin{algorithm}[h!]
\footnotesize
\caption{Decomposition based Branch-and-Bound (DBB)}\label{al:DBB}
\begin{algorithmic}\vspace{0.2cm}
\State \textbf{Step 1}: Initialize $U=\infty$, $k=0$, $\mathcal{P}^k$, and $\mathcal{L}$; 
\\\textbf{Step 2}: \textbf{If} ($\mathcal{L} = \emptyset$) \textbf{then} terminate with solution $\boldsymbol{\hat{\chi}}^*$;
\State \hspace{0.2in} \textbf{Else} select and remove a subproblem $k$ from $\mathcal{L}$ (\textit{i.e.}, $\mathcal{L}=\mathcal{L} \setminus \{k\}$);
\State \textbf{Step 3}: Generate upper ($\gamma^k$) and lower ($\beta^k$) bounds for subproblem $k$;
\State \textbf{Step 4}: Set $L = \text{min}_{i \in \mathcal{L} \cup \{k\}} \beta^i$;
\State \textbf{Step 5}: \textbf{If} ($\gamma^k < U$) \textbf{then} set $U=\alpha^k$ and $\chi^*=\chi^k$;
\State \textbf{Step 6}: Fathom the subproblem (\textit{i.e.}, set $\mathcal{L}=\mathcal{L} \setminus \{i|\beta^i>U\}$);
\State \textbf{Step 7}: \textbf{If} ($\beta^k > U$) \textbf{then} go to Step 2;
\State \textbf{Step 8}: Branch by partitioning $\mathcal{P}^k$ into $\mathcal{P}^{k_1}$ and $\mathcal{P}^{k_2}$;
\State \textbf{Step 9}: Set $\mathcal{L}=\mathcal{L} \cup \{k_1,k_2\}$, $\beta^{k_1}=\beta^k$, $\beta^{k_2}=\beta^k$, $k=k+1$, and go to Step 2.
\end{algorithmic}
\end{algorithm}
In Step 1 of Algorithm \ref{al:DBB}, we begin (after setting $k=0$) by constructing the hyper-rectangle $\mathcal{P}^k = \prod_j(l^k_j,u^k_j] \supset X$ and adding the problem inf\{$f(\chi)|\chi \in X \cap \mathcal{P}^k$\} to the list of open subproblems $\mathcal{L}$. For each component $j$ of $\chi$, we set $l_j=\text{min} \{\chi_j|\chi \in X \}$ and $u_j=\text{max} \{\chi_j|\chi \in X \}$ where the optimization problems are linear programs since $X$ is polyhedral. Then, for each $j$ and scenario $s$, we find $k^s_j \in \mathbb{Z}$ such that $k^s_j-h^s_j-1<l_j<k^s_j-h^s_j$. If $l_j+h^s_j$ is integral then set $k^s_j=l_j+h^s_j$; otherwise $k^s_j=\lfloor l_j+h^s_j+1 \rfloor$. Finally, we set $l^k_j=\text{max}_s\{k^s_j-h^s_j-1 \}$ and $u^k_j=u_j$. In Step 2, to ensure convergence of the algorithm (proof given by \cite{ahmed2004finite}), we select subproblem $k$ such that $\beta^k=L$. In Step 3, for a given subset $\mathcal{P}^k$, we obtain a lower bound on the corresponding subproblem by solving the following formulation:

{\footnotesize
\begin{mini!}|s|[0]
{}
{\theta \label{form:LB}}
{\label{LBform}}
{\text{(LB)} \quad \beta^k=}
\addConstraint{(\ref{TP:c1})}{}{}
\addConstraint{l^k \leq \chi \leq u^k}{}{}
\addConstraint{\theta \geq \sum_{s=1}^{S}{\frac 1 S \Psi^s(l^k + \epsilon)}.}{}{}
\end{mini!}
}
In problem (LB), $\Psi^s(\cdot)$ is constant over $(l^k,l^k+\epsilon]$ for all $s$ when $\epsilon$ is sufficiently small \citep{ahmed2004finite}. The value of $\epsilon$ can be calculated \textit{a priori} using the following algorithm:
\begin{algorithm}[h!]
\footnotesize
\caption{Calculation of $\epsilon$}\label{al:eps}
\begin{algorithmic}\vspace{0.2cm}
\State \textbf{Step 1}: \textbf{For}  each component $j$ of $\chi$
    \State \hspace{0.2in} Set $s=1$, $\Xi=\emptyset$. Choose $k^1_j \in \mathbb{Z}$. Let $\chi^0_j=k^1_j-h^1_j-1$ and $\chi^1_j=\chi^0_j+1$. Set $\Xi=\Xi \cup \{\chi^0_j,\chi^1_j\}$;
\State \textbf{Step 2}: \textbf{For}  $s=2,\ldots,S$
    \State \hspace{0.2in} Set $k^s_j=\lfloor \chi^1_j+h^s_j \rfloor$. Let $\chi^s_j=k^s_j-h^s_j$;
    \State \hspace{0.2in} \textbf{If} $\Xi \cap \{\chi^s_j\}=\emptyset$ \textbf{then} set $\Xi=\Xi \cup \{\chi^s_j\}$;
\State \textbf{Step 3}: Sort the elements of $\Xi$ such that $\chi^0_j=\xi^0_j<\xi^1_j<\ldots<\xi^n_j=\chi^1_j$ with $n \leq S$;
\State \textbf{Step 4}: Let $\epsilon_j=\text{min}_i \{\xi^i_j-\xi^{i-1}_j\}$;
\State \textbf{Step 5}: Set $\epsilon=\frac{1}{2}\text{min}_j \{\epsilon_j\}$.
\end{algorithmic}
\end{algorithm}

In Step 3 of Algorithm \ref{al:DBB}, we also generate an upper bound. For a given subset $\mathcal{P}^k$ such that $\mathcal{P}^k\cap X \neq \emptyset$, let $\chi^k$ be an optimal solution to (LB). Since $\chi^k$ is feasible to (TP) we can simply set $\gamma^k=f(\chi^k)$. Finally, in Step 8 we perform branching. To do this we identify the variable $j^{\prime}$ by determining the value of $\chi_{j^{\prime}}$ where the the current second-stage problem becomes infeasible. For each scenario $s$, let $y^s$ be the solution of the second-stage subproblem when solving (LB). Then, for each $j$ compute $p_j=\text{min}_s \{(D^sy^s)_j-h^s_j\}$. Let $j^\prime \in \text{argmax}_j \{ \min \{ p_j-l^k_j,u^k_j-p_j\}\}$ and split $\mathcal{P}^k$ into $\mathcal{P}^{k_1}=(l^k_{j^\prime},p_{j^\prime}] \prod_{j \neq j^\prime}(l^k_j,u^k_j]$ and $\mathcal{P}^{k_2}=(p_{j^\prime},u^k_{j^\prime}] \prod_{j \neq j^\prime}(l^k_j,u^k_j]$.

\section{Numerical Study}
\label{section:results}
To evaluate the advantages of postponement in making acceptance and scheduling decisions about outpatient appointments, we conducted an extensive numerical study. We also performed a sensitivity analysis to demonstrate how the performance is affected by changes in some problem parameters.
\subsection{Input data}
The patient arrival rates $\lambda_j, \lambda^I, \lambda^E$ and the parameters in Table \ref{table:parameters} are required input for the proposed model. We consider two priority classes of outpatients $(j=1,2)$. The values of $\lambda_j, \lambda^I, \lambda^E$ are estimated based on the average arrival rates of different patient types at the Radiology Department of Prisma Health. The parameters in Table \ref{table:parameters} are estimated with the assistance of physicians at the Radiology Department. However, due to confidentiality concerns we only present normalized values in the appendix. The daily regular capacity of the clinic is estimated to be $K=175$ appointments, and the planning horizon is set to $T=50$ days. In our analysis we ignored the first week (\textit{i.e.}, we used it as our warm-up period). We also ignored the last week of the planning horizon to eliminate any end-of-horizon effects.

\subsection{Experimental setup} 
The problem is implemented in C$^{++}$. The decomposed problems are solved on an Intel Core i7-9700 CPU utilizing the Gurobi 7.0 solver. The computational time required to implement the SAA algorithm will grow as $S$, $S'$, and $M$ increase in Algorithm \ref{al:SAA}. The growth can be linear or exponential depending on whether or not a decomposition approach is used \citep{kleywegt2002sample}. In our case the growth is linear since we are using DBB as presented in Algorithm \ref{al:DBB}. Our first set of experiments were conducted to determine suitable values for $S$, $S'$, and $M$. We began with initializing $S$=10 and $S^{'}$=100, and increased these values in increments of 10 until $\Delta_c$ and $\sigma^2$ values were less than $\epsilon=0.01$. We also tested different values for $M$ from the set $\{10, 20, 30, 40, 50\}$. The final values for $S$, $S'$, and $M$ were respectively, 100, 500, and 30. 

\subsection{The base scenario} \label{section:setup}
After determining the values for $S$, $S'$, and $M$, we conducted a large number of experiments to test and compare the performance of four different appointment scheduling policies. 

\subsubsection{Policy 1:} \label{policy1}
This policy refers to what is currently being used by the clinic. As described earlier, the clinic currently allocates a portion of the regular capacity for emergency patients, a portion to inpatients, and uses the remaining capacity for outpatients. Capacity allocations are done based on the $\alpha^I$ and $\alpha^E$ values obtained from the optimization model without postponement (which is explained below in policy 3). In policy 1, appointment decisions are made as soon as a patient arrives on an FCFS basis and the capacities are dedicated. When an emergency patient arrives that patient is seen immediately if there is capacity, otherwise s/he is referred to another clinic. When an inpatient arrives that patient is given the earliest next day appointment as long as there is capacity, otherwise s/he is seen during overtime hours. When an outpatient arrives s/he is given the earliest possible appointment (regardless of type) over the next seven days. If there is no capacity then s/he is referred to another clinic. 

\subsubsection{Policy 2:} \label{policy2}
This is the proposed policy where the outpatients are kept in an acceptance queue up to 72 hours (3 days). The emergency patients are handled the same way as in policy 1, but the appointment decisions for inpatients and outpatients are made at the end of each day. The SAA and DBB approches are used in policy 2 to determine the capacity allocations in stage 1 and appointment decisions in stage 2. Note that policy 2 is anticipative, \textit{i.e.}, the demand for the whole planning horizon is revealed at the beginning of stage 2. One can think of policy 2 as the policy with perfect information and postponement. 

\subsubsection{Policy 3:} \label{policy3}
This policy is similar to policy 1 in that outpatients are not kept in an acceptance queue. On the other hand, policy 3 is similar to policy 2 because it is anticipative and uses the same SAA and DBB approches to make capacity allocation and appointment scheduling decisions. The main difference is that the regular working hours of a day is split into $T'=54$ periods. In other words, in policy 3 appointment decisions are made every 10 minutes, \textit{i.e.}, near real-time. The decisions regarding acceptance and referral of outpatients are taken in each decision epoch $t'$ $(t'=1,2,...,T')$. The original arrival rates are divided by $T'$ and constraint set (\ref{S2:cons5}) is modified to reflect the fact that unused capacity that is allocated for inpatients cannot be used in policy 3. Additionally, constraint sets (\ref{eqn:Q1}) and (\ref{eqn:Qt}), which capture the evolution of the acceptance queue, are removed from the model. Table \ref{table:parametersestimation2} in the Appendix shows the values of the parameters for policy 3. Policy 3 is essentially the same as policy 1 but with perfect information. 

\subsubsection{Policy 4:} \label{policy4}
Based on our observations of the optimal solutions from policy 2, we developed a simple benchmark policy, which is non-anticipative (\textit{i.e.}, does not rely on knowing the demand for the whole planning horizon) and does not keep patients in an acceptance queue. In policy 4 acceptance or referral of all patients are done on arrival but in a way that mimics the decisions made under policy 2. The optimal values for $\alpha^I$ and $\alpha^E$ obtained from policy 2 are used for capacity allocation. In policy 4 the emergency patients and inpatients are  handled the same way as in policy 1. The outpatients, on the other hand, are handled differently. In policy 1 all outpatients are given earliest available appointments on arrival. In policy 4, however, some of the outpatients are referred to other clinics on arrival regardless of available capacity. As will be shown later, 87.37\% of outpatients that request an appointment end up getting one in policy 2. More specifically, the average acceptance rates are 85\% and 92.5\%, respectively, for Type 1 and Type 2 outpatients. Thus, 15\% (7.5\%) of Type 1 (Type 2) outpatients are immediately referred to another clinic on arrival in policy 4. For those outpatients who are not referred to another clinic, Type 2 outpatients are given the earliest available appointment beginning with day two of the planning period. In other words, next day appointments are not given to Type 2 outpatients. For Type 1 outpatients decisions are made based on $b_{jt^ut^a}$ values. Since policy 2 keeps these patients in the acceptance queue for two days, the $b_{1,2,t^a}$ values are sorted in non-decreasing order, and appointments are given based on this order. However, if the cost difference in consecutive days are within 20\% of each other then the later day in the horizon is selected. For the clinic in our case study, this translates to considering days two, five, and seven from time of arrival for possible appointments. We begin with day two and check the remaining capacity. If this remaining capacity is more than 33\% of the total daily outpatient capacity $K(1-\alpha_I-\alpha_E)$ then the Type 1 outpatient is given an appointment on that day with probability 0.75 (based on our observation of policy 2). For a Type 1 outpatient that does not get an appointment on day two the next option (\textit{i.e.}, day five) is considered and the same rules are applied. This process is repeated until the last day of the planning horizon (day seven in this particular example) at which point all remaining Type 1 outpatients are given an appointment on this last day.


Table \ref{table:scenarios} compares policy 1 to policy 2. As can be see from the table, policy 2 significantly improves the expected average cost for the clinic. The improvement ranged from about 40\% to 46\% depending on how long outpatients are allowed to be kept in the acceptance queue. Recall that the clinic is not willing to keep the outpatients in the acceptance queue more than 3 days. Based on our experiments, the lowest total cost was achieved when $T^u=2$.  Thus, in our remaining experiments the value of $T^u$ is fixed at 2.
\begin{table}[H]
\footnotesize
	\centering
	\caption{Cost improvement and capacity allocation for the base scenario}
	\begin{tabular}{c|c|ccc}
		\hline    
      		& \multicolumn{1}{c|}{Policy 1} & \multicolumn{3}{c}{Policy 2} \\ 
       		& {} & {$T^u=1$} & {$T^u=2$} & {$T^u=3$} \\
		\hline
     		Avg. cost improvement & - & 40.7\% & 45.5\% & 39.9\%   \\
    		\hline
		$\alpha_I$ & 13\% & 20\% & 20\% & 20\% \\
    		\hline
		$\alpha_E$ & 34\% & 34\% & 34\% & 35\% \\
		\hline
 	 \end{tabular}
  	\label{table:scenarios}
\end{table}

Table \ref{table:AllResults} summarizes the results for all four policies. The main takeaway is that policy 4, which is very easy to implement, performs really well. Recall that policy 1 is the current policy used at the clinic, policy 3 is the ``optimal" version of policy 1, policy 2 is the one that keeps outpatients in an acceptance queue, and policy 4 is a simple heuristic that we developed which mimics policy 2. As can be seen from Table \ref{table:AllResults}, policies 2 and 4 result in more than 45\% cost improvement compared to policy 1. The small cost difference between policies 1 and 3 (and between policies 2 and 4) suggest that the value of perfect information is minimal. While this may sound counter-intuitive it is expected, because policies 2 and 3 are run using large samples of scenarios and the average is reported. With respect to emergency patients almost none are referred to another clinic under all four policies. With respect to inpatients policies 2 and 4 handle almost all of them during regular hours, but policies 1 and 3 handle about 3.5\% of them during overtime hours. 

The main difference among the four policies is in the way they handle outpatients. In policy 1, about 82\% of out patients are handled during regular hours and the rest are referred to other clinics. Note that this percentage stays almost the same for Type 1 and Type 2 outpatients which makes sense because the current policy functions on an FCFS basis and does not prioritize outpatients. On the other hand, in policy 3 almost all of Type 2 outpatients are seen during regular hours but only 75\% of Type 1 outpatients are seen during regular hours with an overall average of almost 84\% for all outpatients. Recall that policy 3 is the anticipative version of policy 1. Thus, knowing the demand for the whole planning horizon allows policy 3 to prioritize different outpatients. In policy 2 the percentage of outpatients seen during regular hours is more than 87\%, a relatively significant increase over the current system. Policy 2 is able to do this because it is able to utilize the unused capacity allocated to inpatients. For a fair comparison, the outpatients acceptance percentages for policy 2 include those patients that leave the acceptance queue. For example, if 100 Type 1 outpatients arrive then about 5 leave the acceptance queue and of the remaining 95, on average, 85 get appointments and 10 are referred to other clinics. The performance of policy 4 with respect to patient acceptance is very similar to policy 2 since it was designed to mimic policy 2. 

\begin{table}[H]
\footnotesize
\centering
\caption{Summary of all results for all four policies for the base scenario}
\begin{tabular}{|r|r|r|r|r|}
\hline
 & Policy 1 & Policy 2 & Policy 3 & Policy 4 \\ \hline
Avg. cost improvement & - & 45.5\% & 2.5\% & 45.2\% \\
Emergency patient acceptance &  99.6\% &  99.6\% & 99.6\% & 99.6\% \\
Inpatient acceptance &  96.4\% &  99.9\% & 96.5\% & 99.9\% \\
Outpatient acceptance  &  81.7\% &  87.4\% & 83.9\% & 86.5\% \\ 
Type 1 outpatient acceptance  &  81.5\% &  85.0\% & 75.0\% & 84.2\% \\ 
Type 2 outpatient acceptance  &  82.0\% &  92.5\% & 99.9\% & 92.5\% \\ 
Type 1 outpatients leaving the queue  &  - &  5.0\% & - & - \\ 
Type 2 outpatients leaving the queue  &  - &  5.0\% & - & - \\ 
Inpatient capacity not used  &  4.4\% &  0.0\% & 4.5\% & - \\ 
Days in acceptance queue (Type 1)  & - & 1.9 & - & - \\
Days in acceptance queue (Type 2)  & - & 1.0 & - & - \\
Appointment days ahead (Type 1)  & 6.2 & 5.5 & 6.4 & 5.6 \\
Appointment days ahead (Type 2)  & 6.1 & 3.1 & 6.0 & 3.1 \\
$\alpha^I$ & 13\% & 20\% & 13\% & 20\% \\
$\alpha^E$ & 34\% & 34\% & 34\% & 34\% \\
Solution time (sec.)  & 0.10 & 120 & 110 & 0.14 \\
\hline
\end{tabular}
  \label{table:AllResults}
\end{table}

Table \ref{table:AllResults} also shows the percentage of the overall capacity allocated to each patient group. Recall that policy 1 (policy 4) simply uses the $\alpha^I$ and $\alpha^E$ values obtained from policy 3 (policy 2). Both of policies 2 and 3 allocate about 34\% of the capacity to emergency patients. However, policy 2 allocates more capacity to inpatients compared to policy 3. While the increase from 13\% to 20\% may seem unnecessary, it is expected because under policy 2 with postponement the extra capacity allocated for inpatients can be used for outpatients when needed. 

Another interesting observation is related to the indirect waiting times. As seen in Table \ref{table:AllResults}, under the current policy, the indirect waiting times for Type 1 and Type 2 outpatients are 6.2 and 6.1 days, respectively. As expected, policy 1 does not distinguish between the two types of outpatients. In policy 3, since it is anticipative, the indirect waiting times are 6.4 and 6.0 favoring Type 2 outpatients slightly, but the overall average is almost the same as in policy 1. In policy 2 the indirect waiting times are 7.4 (1.9+5.5) and 4.1 days, respectively, for Type 1 and Type 2 outpatients. This shows that policy 2 prioritizes Type 2 outpatients. The average waiting time is decreased by about 2 days for Type 2 outpatients in the expense of about a 1 day increase for Type 1 outpatients. Policy 4 mimics policy 2 but it does not keep an acceptance queue. Thus, policy 4 does very well in reducing indirect waiting times. 

With respect to computational time, policies 1 and 4 are very fast since they are essentially simulating the appointment system using simple rules. In addition, policies 1 and 4 do not compute capacity allocations but use the values obtained from policies 3 and 2, respectively. Thus, the average CPU time per replication is 0.10 and 0.14 seconds, respectively, for policies 1 and 4. Policies 2 and 3 solve complicated optimization problems to optimize capacity allocations and appointment schedules. As such the average CPU time per replication is 120 and 110 seconds for policies 2 and 3, respectively.

\subsection{Sensitivity analysis}
The results presented in Section \ref{section:setup} demonstrate the effectiveness of policy 2 on the base scenario (referred to as experiment 1). While policy 4 is our proposed policy (because it is easy to implement), it is based on policy 2. Thus, in this section, we perform a sensitivity analysis to observe how policy 2 performs under different conditions. For this analysis, the values of the following parameters are changed one at a time: $b_{jt^ut^a}$, $c^{I}$, $c^{E}$, $\lambda^{I}$, and $p_{jt^u}$. 

\subsubsection{Scheduling costs of outpatients:}
\begin{figure}[htbp]
\centering
\resizebox{1.00\textwidth}{!}{
  \begin{subfigure}{8.75cm}
 \begin{tikzpicture}[baseline]
     \captionsetup{type=figure,justification=centering}
    \begin{axis}[%
    xmin=0, xmax=8,%
    ymin=50, ymax=200,%
    xtick={1,2,3,4,5,6,7},xticklabels={1,2,3,4,5,6,7},%
    ytick={50,75,100,125,150,175,200},yticklabels={50,75,100,125,150,175,200},%
     xlabel=Experiment Number,
   ylabel=Total cost,
    ]
    \boxplot{1}{100}{90.472}{110.603}{57.03}{157.065}
     \boxplot{2}{125.237}{114.643}{135.736}{78.97}{184.916}
       \boxplot{3}{132.351}{121.76}{141.482}{86.795}{183.524}
     \boxplot{4}{117.670}{107.063}{128.139}{73.424}{177.658}
       \boxplot{5}{119.985}{109.205}{130.603}{74.374}{181.036}
       \boxplot{6}{120.53}{109.671}{130.508}{73.488}{179.557}
      \boxplot{7}{126.157}{115.553}{136.553}{79.378}{183.03}
    \end{axis}
    \end{tikzpicture}
     \caption{}
          \label{fig:outcost1}
  \end{subfigure}
  \begin{subfigure}{8.75cm}
\centering
\begin{tikzpicture}
\begin{axis}[
   xlabel=Experiment Number,
    ylabel=Capacity allocation,
    xmin=0, xmax=8,
    ymin=0.1, ymax=0.4,
    xtick={1,2,3,4,5,6,7},
    ytick={0,.1,.2,.3,.4,.5},
    legend pos=north east,
    ymajorgrids=true,
    grid style=dashed,
]
\addplot+[line width=2pt,mark size=3pt][
    color=blue,
    mark=triangle,
    ]
    coordinates {
    (1,.2)(2,.2)(3,.22)(4,.2)(5,.2)(6,.2)(7,.2)};
    \addplot+[line width=2pt,mark size=3pt][
    color=red,
    mark=square,
    ]
    coordinates {
    (1,.34)(2,.34)(3,.34)(4,.34)(5,.34)(6,.34)(7,.34)};
    \legend{$\alpha_I$,$\alpha_E$}
	\addlegendentry{estimate}
\end{axis}
\end{tikzpicture}
\caption{}
    \label{fig:outcost2}
  \end{subfigure}
  }
  \caption{Effect of changing $b_{jt^ut^a}$ on average total cost and  capacity allocation}
    \label{fig:outcost}
\end{figure}
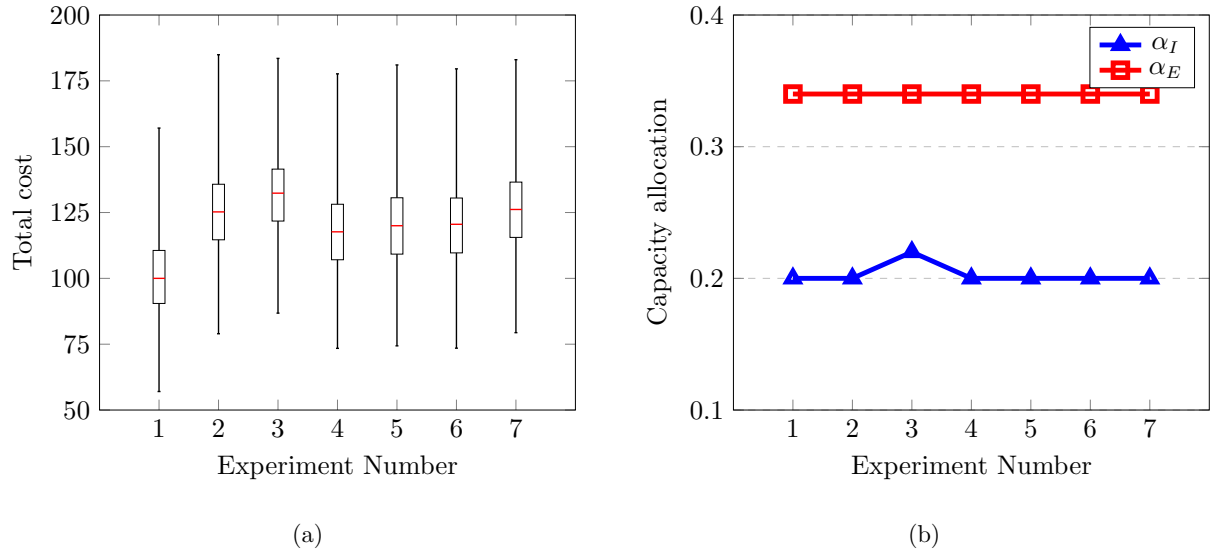

\begin{figure}[htbp]
\centering
\resizebox{1.00\textwidth}{!}{
  \begin{subfigure}{8.75cm}
\centering
\begin{tikzpicture}
\begin{axis}[
   xlabel=Experiment Number,
    ylabel=Days in acceptance queue,
    xmin=0, xmax=6,
    ymin=0, ymax=3,
    xtick={1,2,3,4,5},
    ytick={0,1,2,3},
    legend pos=north east,
    ymajorgrids=true,
    grid style=dashed,
]
\addplot+[line width=2pt,mark size=3pt][
    color=blue,
    mark=triangle,
    ]
    coordinates {
    (1,1.98)(2,1.93)(3,1.83)(4,1.08)(5,1.01)};
    \addplot+[line width=2pt,mark size=3pt][
    color=red,
    mark=square,
    ]
    coordinates {
    (1,1)(2,1)(3,1)(4,1)(5,1)};
    \legend{Type 1,Type 2}
\end{axis}
\end{tikzpicture}
\caption{}
\label{fig:2outcost1}
  \end{subfigure}
  \begin{subfigure}{8.75cm}
\centering
\begin{tikzpicture}
\begin{axis}[
   xlabel=Experiment Number,
    ylabel=Appointment days ahead,
    xmin=0, xmax=6,
    ymin=0, ymax=6,
    xtick={1,2,3,4,5},
    ytick={0,1,2,3,4,5,6},
     xticklabels={1,2,3,6,7},%
    legend pos=north east,
    ymajorgrids=true,
    grid style=dashed,
]
\addplot+[line width=2pt,mark size=3pt][
    color=blue,
    mark=triangle,
    ]
    coordinates {
    (1,5.48)(2,4.6)(3,2.65)(4,3.36)(5,2.38)};
    \addplot+[line width=2pt,mark size=3pt][
    color=red,
    mark=square,
    ]
    coordinates {
    (1,3.08)(2,1.29)(3,1)(4,2.1)(5,2.1)};
    \legend{Type 1,Type 2}
\end{axis}
\end{tikzpicture}
\caption{}
\label{fig:2outcost2}
  \end{subfigure}
  }
  \caption{Effect of changing $b_{jt^ut^a}$ on the number of days in acceptance queue and indirect waiting time}
    \label{fig:2outcost}
\end{figure}
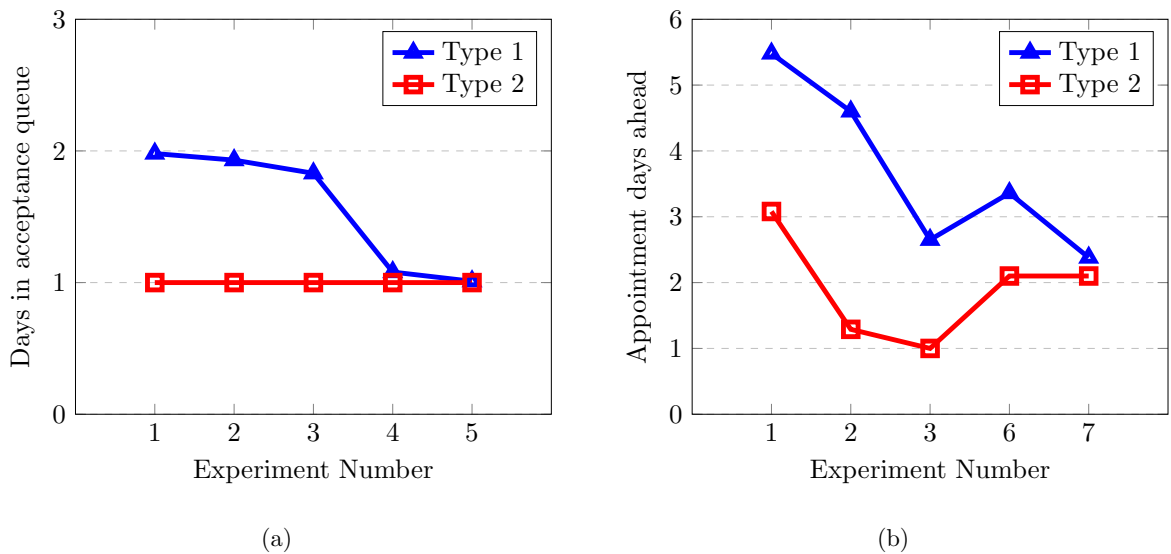

To understand the effect of changing outpatient scheduling costs on the optimal solution, we increased $b_{jt^ut^a}$ by $50\%$ and $100\%$ in experiments $2$ and $3$, respectively. By increasing all the $b_{jt^ut^a}$ values with the same percentage we penalize both wait times in the acceptance queue and the time until appointments in experiments 2 and 3. In experiments $4$ and $5$, we penalize long wait times in the acceptance queue by increasing the $b_{jt^ut^a}$ values for only $t^u=2$ by $50\%$ and $100\%$, respectively. In experiments $6$ and $7$, we penalize scheduling later appointments, where $b_{jt^ut^a}$ values for $t^a \geq 3$ are increased by $50\%$ and $100\%$, respectively. Figures \ref{fig:outcost} and \ref{fig:2outcost} represent the results of experiments 1-7. As seen in Figure \ref{fig:outcost2}, the percent of capacity allocated for emergency patients ($\alpha^E=0.34$) is not impacted by changes to $b_{jt^ut^a}$. On the other hand, capacity allocated for inpatients ($\alpha^I$) increases slightly from 20\% to 22\% in experiment 3 since we are willing to reserve more next day appointments for outpatients. The box plots in Figure \ref{fig:outcost1} show the normalized values of average total cost for the system. In other words, the average total cost for experiment 1 is normalized to 100, and thus, the other values show the corresponding change in cost. As expected, the total cost increases the most in experiments 2 and 3 since all $b_{jt^ut^a}$ values are increased here whereas only a subset of the $b_{jt^ut^a}$ values are increased in experiments 4-7. 

Figure \ref{fig:2outcost1} shows that in the base scenario, all Type 2 outpatients wait in the acceptance queue for one day before they are given an appointment. Type 1 outpatients, however, wait for almost two days in the queue. As the cost of waiting in the queue increases the waiting time for Type 2 outpatients remain the same. For Type 1 outpatients, it decreases. In other words, acceptance and referral decisions are made sooner. Figure \ref{fig:2outcost2} shows that time from the day of acceptance to the day of appointment decreases as $b_{jt^ut^a}$ increases. Because all $b_{jt^ut^a}$ are increased in experiments 2 and 3, the indirect waiting time decreases sharply. In experiments 6 and 7, the $b_{jt^ut^a}$ values were increased for only high values of $t^a$, as such, compared to the base case the drop in indirect waiting time is not as dramatic. However, policy 2 tries to offer earlier appointments to both outpatient types in experiments 6 and 7. Additional insight on these experiments are also discussed in Section \ref{MoreInsight} based on Tables \ref{table:prioritizing} and \ref{table:details}.

\subsubsection{Referral and overtime costs of emergency patients and inpatients: }
Emergency patients are the highest priority patients followed by inpatients. Lack of available capacity to schedule them during regular hours results in additional cost, specifically, overtime cost $c^I$ for inpatients and referral cost $c^E$ for emergency patients. Experiments 8 and 9 measure the effect of increasing these parameters by 50\% and 100\%, respectively. With respect to capacity allocation, $\alpha^E$ and $\alpha^I$ remained the same in experiments 8 and 9 as they were in experiment 1. The only difference was on the average total cost, which increased by about 4\% from experiment 1 to 8 and about 8\% from experiment 1 to 9.       


\subsubsection{Inpatient arrivals:}
Since capacity allocated to inpatients can also be used for outpatients in policy 2, arrival rate of inpatients affect the scheduling of outpatients. To observe this impact, we performed experiments $10$ and $11$ where $\lambda^I$ is increased by $25\%$ and $50\%$, respectively. As one would expect, increasing the inpatient arrival rate resulted in higher $\alpha^I$ values (25\% in experiment 10 and 28\% in experiment 11). On the other hand, there was no change to the capacity allocated to emergency patients which stayed at 34\%. However, the outpatient acceptance rate decreased in both of these experiments. Since more inpatients are arriving into the system, there is less capacity left for outpatients and more of them are referred to other clinics. Given that all problem parameters (including total capacity) remained the same as $\alpha^I$ was increased, the total system cost increased about $28\%$ and $44\%$ in experiments $10$ and $11$, respectively, compared to experiment 1.   

\subsubsection{Abandonment rate of outpatients:}
As discussed earlier, one of the disadvantages of implementing postponed acceptance in a service system such as a diagnostic clinic is that customers may abandon the acceptance queue. We assume that a proportion (1 - $p_{jt^u}$) of type $j$ outpatients leave the queue after having waited $t^u$ days. In addition to the outpatient type and the amount of time they have waited in the queue, the abandonment rates also depend on the type of diagnostic clinic. To capture the effect of such changes, we decreased the value of $p_{jt^u}$ by 5\% in experiment 12 and 10\% in experiment 13. In other words, the chances of an outpatient abandoning the acceptance queue is higher in experiments 12 and 13. 
As a result, $\alpha^I$ stayed at 20\% in experiment 12 but slightly increased to 21\% in experiment 13. Since fewer outpatients are waiting for an appointment due to abandonment, the rejection rate of outpatients decreased in experiments 12 and 13 compared to experiment 1. More specifically, in experiment 12 the number of outpatients referred to other clinics decreased by $85\%$. In experiment 13 no outpatient was referred to another clinic. This resulted in lower average total cost with a decrease of about $27\%$ and $40\%$ in experiments 12 and 13, respectively, as compared to experiment 1.

\subsubsection{Additional insights:} \label{MoreInsight}
Table \ref{table:prioritizing} provides the percentage of outpatients who have received appointments after waiting one day or two days in the acceptance queue. Type 2 outpatients always received an appointment after only one day in the acceptance queue. 
On the other hand, majority of type 1 outpatients wait for two days in the acceptance queue in all of the experiments except experiments 4 and 5. Recall that in experiments 4 and 5 the $b_{jt^ut^a}$ values are increased for $t^u=2$ by $50\%$ and $100\%$, respectively. In other words, waiting in the acceptance queue for two days is costly in these cases. Thus, in experiments 4 and 5 most of the outpatients get an appointment after only one day in the queue. 

\begin{table}[h]
\footnotesize
\centering
\caption{Percentage of outpatients waiting in the acceptance queue for one vs. two days}
\label{table:prioritizing}
\begin{tabular}{|c|c|c|c|c|}
\hline
  \begin{tabular}[c]{@{}c@{}} \\Experiment\end{tabular} & \begin{tabular}[c]{@{}c@{}}Type 1 outpatients \\One day in queue\end{tabular} & \begin{tabular}[c]{@{}c@{}}Type 1 outpatients \\ Two days in queue\end{tabular} & \begin{tabular}[c]{@{}c@{}}Type 2 outpatients \\ One day in queue\end{tabular} & \begin{tabular}[c]{@{}c@{}}Type 2 outpatients \\ Two days in queue\end{tabular}\\ \hline
 1 & 1.93 & 98.07 & 100.00 & 0.00 \\
 2 & 6.63& 93.37 &100.00 & 0.00 \\
 3  & 16.60 & 83.40 &100.00 & 0.00 \\
 4  & 92.14  & 7.86 & 100.00 & 0.00\\
  5  &  99.40 & 0.60 & 100.00 & 0.00\\
  6  & 7.61  & 92.39 & 100.00 & 0.00\\
  7  & 1.97 & 98.03 & 100.00 & 0.00\\
  8  & 1.70  & 98.44 & 100.00 & 0.00\\
  9  & 1.56 & 98.43 & 100.00 & 0.00\\
  10  & 0.54 & 99.46 & 100.00 & 0.00\\
  11  & 0.22 & 98.78 & 100.00 & 0.00\\
   12  & 2.18 & 97.82 & 100.00 & 0.00\\
  13  & 3.59 & 96.41& 100.00 & 0.00\\
  \hline
\end{tabular}
\end{table}

The unused inpatient capacity before and after scheduling outpatients out of the acceptance queue are captured and listed in columns 2 and 3 of Table \ref{table:details}, respectively. Columns 4, 5, and 6 of the table shows the percentage of the patients that ultimately received appointments. Depending on the problem parameters, the unused inpatient capacity varies between 44\% and 56\%. Note that almost all of the inpatients receive appointments and the leftover capacity is used for outpatients. Outpatient acceptance rate is low for experiment 11 but over 80\% in all of the other experiments. Recall that in experiments 10 and 11 the inpatient arrival rate was increased, as such, there is not much leftover capacity that can be used for outpatients compared to the other experiments. 

\begin{table}[h]
\footnotesize
\centering
\caption{Capacity utilization of each patient type for policy 2}
\label{table:details}
\resizebox{\textwidth}{!}{
\begin{tabular}{|c|c|c|c|c|c|}
\hline
  \begin{tabular}[c]{@{}c@{}}Experiment\end{tabular} &  \begin{tabular}[c]{@{}c@{}}Unused inpatient capacity \\ before scheduling \\ outpatients (\%)\end{tabular} &\begin{tabular}[c]{@{}c@{}}Unused inpatient capacity \\ after scheduling \\ outpatients (\%)\end{tabular}  & \begin{tabular}[c]{@{}c@{}} Inpatient \\ acceptance (\%)\end{tabular} & \begin{tabular}[c]{@{}c@{}} Emergency patient \\ acceptance (\%)\end{tabular} & \begin{tabular}[c]{@{}c@{}} Outpatient \\ acceptance (\%)\end{tabular}  \\ \hline
 1 &  44.45 & 0.00 & 99.99 &  99.56 & 87.39 \\
 2 &  44.45 & 0.00 & 99.99 & 99.56  & 85.15 \\
 3 & 52.78 & 0.00 & 99.99 & 99.56 & 82.90 \\
 4 & 44.45 & 0.00 & 99.99 & 99.56 & 86.08  \\
  5 & 44.45 & 0.00 & 99.99  & 99.56 & 84.41 \\
  6 & 44.45 & 0.00 & 99.99  & 99.56  & 85.25 \\
  7 & 44.45 & 0.00 & 99.99  & 99.56  & 82.38 \\
  8 & 44.45 & 0.00 & 99.99  & 99.56  & 87.40  \\
  9 & 44.45 & 0.00 &  99.99 & 99.56  & 87.41\\
  10 & 52.75 & 0.00 & 99.99  & 99.55  & 81.49 \\
  11 & 55.47 & 0.00 & 99.99 & 99.57  & 75.74  \\
   12 & 44.45 & 0.00 & 99.99 & 99.56 & 84.31 \\
  13 & 50.00 & 0.00 & 99.99 & 99.56  & 84.03 \\
  \hline
\end{tabular}
}
\end{table}

\section{Conclusion} \label{section:conclusion}
This paper introduces a postponable acceptance appointment system for a diagnostic clinic. Diagnostic facilities often serve patients of different priority classes. Outpatients are typically scheduled in advance, but higher priority patients (\textit{i.e.}, inpatients and emergency patients) are usually seen as soon as possible. Scheduling of outpatients at the clinic are currently  done on an FCFS basis. Thus, high priority outpatients may not receive timely appointments. This challenge motivated us to propose a postponement system in scheduling of different patient classes. The value of the proposed model is that the system can strategically postpone the acceptance of low priority outpatients while waiting for higher priority outpatients. We formulate the problem as a TSSP model in which the first stage estimates the optimal capacity reserved for inpatients and emergency patients. In the second stage, the decisions regarding acceptance and referral of outpatients are made. 

Using a data set from the Radiology Department of Prisma Health, we have conducted a series of experiments to test how the model works.
The results suggest that postponing the acceptance or referral of outpatient appointment requests up to two days improves the system-wide cost while reducing indirect waiting times. The cost improvement achieved is primarily due to the increase in the utilization of the unused inpatient capacity for outpatients waiting in the queue. In addition, the system prioritizes more urgent outpatients by having them wait only one day in the queue and forcing the less urgent outpatients to wait for two days in the acceptance queue. After analyzing the optimal solutions obtained from our model we developed a simple benchmark policy that can be implemented in real life which performs well.

This study can be extended in multiple directions. For example, in this study we assume that the duration of visits are constant and identical for each type of patient. Thus, the number of patients that can be seen each day is a fixed number. To consider a more realistic case, uncertain service times can be considered. Furthermore, due to the higher indirect waiting time of lower priority patients, the possibility of no-shows may increase for these classes. Thus, the model can be extended to consider no-shows. Another extension could be to develop a multi-stage stochastic programming approach since the demand uncertainty is revealed over time after each time period. Such multi-stage approaches will be computationally more difficult to solve. Alternatively, the two-stage stochastic program can be used on a rolling horizon basis. 
\newpage
\section{Appendix}
\begin{table}[h]
\footnotesize
\centering
 \caption{Parameters values for policy 2}
  \label{table:parametersestimation}
\begin{tabular}{|c|c||c|c||c|c|}
\hline
 Parameter& Value & Parameter & Value & Parameter & Value  \\ \hline
$\lambda^E$ & 50 & $b_{1,2,4}$ & 0.75 &  $b_{2,2,5}$ & 2 \\ 
$\lambda^I$ & 20 & $b_{1,2,5}$ & 0.75 & $b_{2,2,6}$ & 2.5 \\
$\lambda_1$ & 80 & $b_{1,2,6}$ & 1.25 & $b_{2,2,7}$ & 2.5 \\ 
$\lambda_2$ & 40 & $b_{1,2,7}$ & 1.25 & $b_{2,3,1}$ & 3 \\ 
$c^E$ & 45 & $b_{1,3,1}$ & 1.5 & $b_{2,3,2}$ & 4 \\ 
$c^I$& 30 & $b_{1,3,2}$ & 2.5 & $b_{2,3,3}$ & 4 \\ 
$c^O_{1,1}$ & 8 & $b_{1,3,3}$ & 2.5 & $b_{2,3,4}$ &5 \\ 
$c^O_{1,2}$ & 12 & $b_{1,3,4}$ & 3.5 & $b_{2,3,5}$ & 5 \\ 
$c^O_{1,3}$ & 16 & $b_{1,3,5}$ & 3.5 & $b_{2,3,6}$ & 6 \\  
$c^O_{2,1}$ & 14 & $b_{1,3,6}$ & 4.5 & $b_{2,3,7}$ & 6 \\ 
$c^O_{2,2}$ & 18 & $b_{1,3,7}$ & 4.5 & $a_{1,1}$ & 1 \\ 
$c^O_{2,3}$ & 22 & $b_{2,1,1}$ & 0 & $a_{1,2}$ & 1.25\\ 
$b_{1,1,1}$ & 0 & $b_{2,1,2}$ & 0.5 & $a_{1,3}$ & 4.5 \\ 
$b_{1,1,2}$ & 0 & $b_{2,1,3}$ & 0.5 & $a_{2,1}$ & 1.5\\ 
$b_{1,1,3}$ & 0.5 & $b_{2,1,4}$ & 1 & $a_{2,2}$ & 2.5\\ 
$b_{1,1,4}$ & 0.5 & $b_{2,1,5}$ & 1 & $a_{2,3}$ & 6\\ 
$b_{1,1,5}$ & 0.5 & $b_{2,1,6}$ & 1.5 & $p_{1,1}$ & 1\\ 
$b_{1,1,6}$ & 1 & $b_{2,1,7}$& 1.5 & $p_{1,2}$ & 0.95\\ 
$b_{1,1,7}$ & 1 & $b_{2,2,1}$ & 1 & $p_{1,3}$ & 0.9\\ 
$b_{1,2,1}$ & 0.25 & $b_{2,2,2}$ & 1.5 & $p_{2,1}$ & 0.95\\ 
$b_{1,2,2}$ & 0.25 & $b_{2,2,3}$ & 1.5 & $p_{2,2}$ & 0.9\\ 
$b_{1,2,3}$ & 0.75 & $b_{2,2,4}$ & 2 & $p_{2,3}$ & 0.85 \\ 
\hline
\end{tabular}
\end{table}

\begin{table}[h]
\footnotesize
\centering
\caption{Parameter values for policy 3}
  \label{table:parametersestimation2}
\begin{tabular}{|c|c||c|c||c|c|}
\hline
Parameter& Estimation & Parameter & Estimation & Parameter & Estimation  \\ \hline
$\lambda^E$ & 50 & $b_{1,1}$ & 0 & $b_{2,2}$ & 0.25  \\ 
$\lambda^I$ & 20 & $b_{1,2}$ & 0 & $b_{2,3}$ & 0.25  \\
$\lambda_1$ & $\frac{80}{54}$ & $b_{1,3}$ & 0.25 & $b_{2,4}$ & 0.5 \\ 
$\lambda_2$ & $\frac{40}{54}$ & $b_{1,4}$ & 0.25 & $b_{2,5}$ & 0.5  \\ 
$c^E$ & 45 & $b_{1,5}$ & 0.25 & $b_{2,6}$ & 0.75 \\ 
$c^I$& 30 & $b_{1,6}$ & 0.5 & $b_{2,7}$ &0.75 \\ 
$c^O_{1}$ & 6 & $b_{1,7}$ & 0.5 & & \\ 
$c^O_{2}$ & 10 & $b_{2,1}$ & 0 & & \\  
\hline
\end{tabular}
\end{table}

\newpage
\bibliographystyle{plain}
\bibliography{main}

\begin{thebibliography}{35}
\providecommand{\natexlab}[1]{#1}
\providecommand{\url}[1]{\texttt{#1}}
\providecommand{\urlprefix}{URL }

\bibitem[{Ahmadi-Javid et~al.(2017)Ahmadi-Javid, Jalali, \protect\BIBand{}
  Klassen}]{ahmadi}
Ahmadi-Javid A, Jalali Z, Klassen KJ (2017) Outpatient appointment systems in
  healthcare: A review of optimization studies. \emph{European Journal of
  Operational Research} 258(1):3--34.

\bibitem[{Ahmed et~al.(2002)Ahmed, Shapiro, \protect\BIBand{}
  Shapiro}]{ahmed2002sample}
Ahmed S, Shapiro A, Shapiro E (2002) The sample average approximation method
  for stochastic programs with integer recourse. \emph{Submitted for
  publication} 1--24.

\bibitem[{Ahmed et~al.(2004)Ahmed, Tawarmalani, \protect\BIBand{}
  Sahinidis}]{ahmed2004finite}
Ahmed S, Tawarmalani M, Sahinidis NV (2004) A finite branch-and-bound algorithm
  for two-stage stochastic integer programs. \emph{Mathematical Programming}
  100(2):355--377.

\bibitem[{Akhavizadegan et~al.(2017)Akhavizadegan, Ansarifar, \protect\BIBand{}
  Jolai}]{akhavizadeganl}
Akhavizadegan F, Ansarifar J, Jolai F (2017) A novel approach to determine a
  tactical and operational decision for dynamic appointment scheduling at
  nuclear medical center. \emph{Computers \& Operations Research} 78:267--277.

\bibitem[{Balasubramanian et~al.(2013)Balasubramanian, Muriel, Ozen, Wang, Gao,
  \protect\BIBand{} Hippchen}]{balasubramanian}
Balasubramanian H, Muriel A, Ozen A, Wang L, Gao X, Hippchen J (2013) Capacity
  allocation and flexibility in primary care. \emph{Handbook of healthcare
  operations management}, 205--228 (Springer).

\bibitem[{Bassamboo et~al.(2005)Bassamboo, Harrison, \protect\BIBand{}
  Zeevi}]{bassamboo}
Bassamboo A, Harrison JM, Zeevi A (2005) Dynamic routing and admission control
  in high-volume service systems: Asymptotic analysis via multi-scale fluid
  limits. \emph{Queueing Systems} 51(3-4):249--285.

\bibitem[{Berg et~al.(2014)Berg, Denton, Erdogan, Rohleder, \protect\BIBand{}
  Huschka}]{berg2014optimal}
Berg BP, Denton BT, Erdogan SA, Rohleder T, Huschka T (2014) Optimal booking
  and scheduling in outpatient procedure centers. \emph{Computers \& Operations
  Research} 50:24--37.

\bibitem[{Bertsimas \protect\BIBand{} Popescu(2003)}]{bertsimas2003revenue}
Bertsimas D, Popescu I (2003) Revenue management in a dynamic network
  environment. \emph{Transportation science} 37(3):257--277.

\bibitem[{Bhattacharjee \protect\BIBand{} Ray(2016)}]{bhattacharjee}
Bhattacharjee P, Ray PK (2016) Simulation modelling and analysis of appointment
  system performance for multiple classes of patients in a hospital: a case
  study. \emph{Operations Research for Health Care} 8:71--84.

\bibitem[{Deglise-Hawkinson et~al.(2018)Deglise-Hawkinson, Helm, Huschka,
  Kaufman, \protect\BIBand{} Van~Oyen}]{deglise2018capacity}
Deglise-Hawkinson J, Helm JE, Huschka T, Kaufman DL, Van~Oyen MP (2018) A
  capacity allocation planning model for integrated care and access management.
  \emph{Production and operations management} 27(12):2270--2290.

\bibitem[{Erdogan \protect\BIBand{} Denton(2013)}]{erdogan2013dynamic}
Erdogan SA, Denton B (2013) Dynamic appointment scheduling of a stochastic
  server with uncertain demand. \emph{INFORMS Journal on Computing}
  25(1):116--132.

\bibitem[{Erdogan et~al.(2015)Erdogan, Gose, \protect\BIBand{}
  Denton}]{erdogan2015online}
Erdogan SA, Gose A, Denton BT (2015) Online appointment sequencing and
  scheduling. \emph{IIE Transactions} 47(11):1267--1286.

\bibitem[{Feldman et~al.(2014)Feldman, Liu, Topaloglu, \protect\BIBand{}
  Ziya}]{feldman2014appointment}
Feldman J, Liu N, Topaloglu H, Ziya S (2014) Appointment scheduling under
  patient preference and no-show behavior. \emph{Operations Research}
  62(4):794--811.

\bibitem[{Gao et~al.(2012)Gao, Xu, \protect\BIBand{} Ball}]{gao2012managing}
Gao L, Xu SH, Ball MO (2012) Managing an available-to-promise assembly system
  with dynamic short-term pseudo-order forecast. \emph{Management Science}
  58(4):770--790.

\bibitem[{Green et~al.(2006)Green, Savin, \protect\BIBand{} Wang}]{green}
Green LV, Savin S, Wang B (2006) Managing patient service in a diagnostic
  medical facility. \emph{Operations Research} 54(1):11--25.

\bibitem[{Gupta \protect\BIBand{} Denton(2008)}]{gupta}
Gupta D, Denton B (2008) Appointment scheduling in health care: Challenges and
  opportunities. \emph{IIE transactions} 40(9):800--819.

\bibitem[{Jiang et~al.(2017)Jiang, Shen, \protect\BIBand{}
  Zhang}]{jiang2017integer}
Jiang R, Shen S, Zhang Y (2017) Integer programming approaches for appointment
  scheduling with random no-shows and service durations. \emph{Operations
  Research} 65(6):1638--1656.

\bibitem[{Kang et~al.(2016)Kang, Shanthikumar, \protect\BIBand{}
  Altinkemer}]{kang}
Kang K, Shanthikumar JG, Altinkemer K (2016) Postponable acceptance and
  assignment: A stochastic dynamic programming approach. \emph{Manufacturing \&
  Service Operations Management} 18(4):493--508.

\bibitem[{Kesling \protect\BIBand{} Nissenbaum(2014)}]{kesling}
Kesling B, Nissenbaum D (2014) {VA Goal to Slash Wait Times Was
  ‘Unrealistic’, Aide Said}. \emph{The Wall Street Journal (May 23).
  Available at https://www. wsj.
  com/articles/SB10001424052702303749904579580473122138420 (accessed date
  August 27, 2017)} .

\bibitem[{Kleywegt et~al.(2002)Kleywegt, Shapiro, \protect\BIBand{} Homem-de
  Mello}]{kleywegt2002sample}
Kleywegt AJ, Shapiro A, Homem-de Mello T (2002) The sample average
  approximation method for stochastic discrete optimization. \emph{SIAM Journal
  on Optimization} 12(2):479--502.

\bibitem[{Kong et~al.(2015)Kong, Li, Liu, Teo, \protect\BIBand{}
  Yan}]{kong2015appointment}
Kong Q, Li S, Liu N, Teo CP, Yan Z (2015) Appointment scheduling under
  schedule-dependent patient no-show behavior.

\bibitem[{Laporte \protect\BIBand{} Louveaux(1993)}]{laporte1993integer}
Laporte G, Louveaux FV (1993) The integer l-shaped method for stochastic
  integer programs with complete recourse. \emph{Operations research letters}
  13(3):133--142.

\bibitem[{Luo et~al.(2012)Luo, Kulkarni, \protect\BIBand{} Ziya}]{luo}
Luo J, Kulkarni VG, Ziya S (2012) Appointment scheduling under patient no-shows
  and service interruptions. \emph{Manufacturing \& Service Operations
  Management} 14(4):670--684.

\bibitem[{Mak et~al.(1999)Mak, Morton, \protect\BIBand{} Wood}]{mak}
Mak WK, Morton DP, Wood RK (1999) Monte carlo bounding techniques for
  determining solution quality in stochastic programs. \emph{Operations
  research letters} 24(1-2):47--56.

\bibitem[{McCarthy et~al.(2000)McCarthy, McGee, \protect\BIBand{}
  O'Boyle}]{mccarthy}
McCarthy K, McGee H, O'Boyle C (2000) Outpatient clinic waiting times and
  non-attendance as indicators of quality. \emph{Psychology, health \&
  medicine} 5(3):287--293.

\bibitem[{Moore et~al.(2001)Moore, Wilson-Witherspoon, \protect\BIBand{}
  Probst}]{moore}
Moore CG, Wilson-Witherspoon P, Probst JC (2001) Time and money: effects of
  no-shows at a family practice residency clinic. \emph{Family Medicine-Kansas
  City-} 33(7):522--527.

\bibitem[{Patrick \protect\BIBand{} Puterman(2007)}]{patrick2007improving}
Patrick J, Puterman ML (2007) Improving resource utilization for diagnostic
  services through flexible inpatient scheduling: A method for improving
  resource utilization. \emph{Journal of the Operational Research Society}
  58(2):235--245.

\bibitem[{Patrick et~al.(2008)Patrick, Puterman, \protect\BIBand{}
  Queyranne}]{patrick}
Patrick J, Puterman ML, Queyranne M (2008) Dynamic multipriority patient
  scheduling for a diagnostic resource. \emph{Operations research}
  56(6):1507--1525.

\bibitem[{Qu et~al.(2013)Qu, Peng, Kong, \protect\BIBand{} Shi}]{qu2013two}
Qu X, Peng Y, Kong N, Shi J (2013) A two-phase approach to scheduling
  multi-category outpatient appointments--a case study of a women’s clinic.
  \emph{Health care management science} 16(3):197--216.

\bibitem[{Sickinger \protect\BIBand{} Kolisch(2009)}]{sickinger}
Sickinger S, Kolisch R (2009) The performance of a generalized bailey--welch
  rule for outpatient appointment scheduling under inpatient and emergency
  demand. \emph{Health care management science} 12(4):408.

\bibitem[{Smith et~al.(1992)Smith, Leimkuhler, \protect\BIBand{}
  Darrow}]{smith1992yield}
Smith BC, Leimkuhler JF, Darrow RM (1992) Yield management at american
  airlines. \emph{interfaces} 22(1):8--31.

\bibitem[{Tang \protect\BIBand{} Wang(2015)}]{tang}
Tang J, Wang Y (2015) An adjustable robust optimisation method for elective and
  emergency surgery capacity allocation with demand uncertainty.
  \emph{International Journal of Production Research} 53(24):7317--7328.

\bibitem[{Van~Slyke \protect\BIBand{} Wets(1969)}]{van1969shaped}
Van~Slyke RM, Wets R (1969) L-shaped linear programs with applications to
  optimal control and stochastic programming. \emph{SIAM Journal on Applied
  Mathematics} 17(4):638--663.

\bibitem[{Verweij et~al.(2003)Verweij, Ahmed, Kleywegt, Nemhauser,
  \protect\BIBand{} Shapiro}]{verweij}
Verweij B, Ahmed S, Kleywegt AJ, Nemhauser G, Shapiro A (2003) The sample
  average approximation method applied to stochastic routing problems: a
  computational study. \emph{Computational Optimization and Applications}
  24(2-3):289--333.

\bibitem[{Weatherford \protect\BIBand{} Bodily(1992)}]{weatherford1992taxonomy}
Weatherford LR, Bodily SE (1992) A taxonomy and research overview of
  perishable-asset revenue management: Yield management, overbooking, and
  pricing. \emph{Operations research} 40(5):831--844.

\end{thebibliography}

\end{document}